\documentclass[aos,MSNbibl,nameyear,rotating,seceqn,dvips]{arximspdf}
\usepackage{dcolumn}

\doi{10.1214/11-AOS897}
\volume{39}
\issue{4}
\pubyear{2011}
\firstpage{2021}
\lastpage{2046}

\makeatletter

\def\bptnote#1{}

\newcolumntype{d}[1]{D{.}{.}{#1}}

\newcommand{\ba}{\mathbf{a}}
\newcommand{\bb}{\mathbf{b}}
\newcommand{\hb}{\hat{b}}
\newcommand{\cE}{\mathcal{E}}
\newcommand{\cG}{\mathcal{G}}
\newcommand{\bI}{\mathbf{I}}
\newcommand{\tM}{\widetilde{M}}
\newcommand{\cO}{\mathcal{O}}
\newcommand{\rP}{\mathrm{P}}
\newcommand{\tP}{\widetilde{P}}
\newcommand{\tbr}{\tilde{\mathbf{r}}}
\newcommand{\bu}{\mathbf{u}}
\newcommand{\bv}{\mathbf{v}}
\newcommand{\bx}{\mathbf{x}}
\newcommand{\tX}{\widetilde{X}}
\newcommand{\by}{\mathbf{y}}
\newcommand{\tby}{\tilde{\mathbf{y}}}
\newcommand{\bzero}{\mathbf{0}}
\newcommand{\eps}{\epsilon}
\newcommand{\veps}{\varepsilon}

\newcommand{\argmin}{\mathop{\arg\min}}
\newcommand{\real}{\mathbb{R}}
\newcommand{\diag}{\operatorname{diag}}

\newcommand{\sgn}{\operatorname{sgn}}

\newcommand{\hbeta}{\hat{\beta}}
\newcommand{\trace}{\operatorname{trace}}

\newcommand{\bone}{\mathbf{1}}

\newcommand{\bbeta}{\bolds\beta}
\newcommand{\hbbeta}{\hat{\bolds\beta}}
\newcommand{\tbeta}{\tilde{\beta}}
\newcommand{\tbbeta}{\tilde{\bolds\beta}}
\newcommand{\bdelta}{\bolds\delta}
\newcommand{\gam}{\gamma}
\newcommand{\lam}{\lambda}
\newcommand{\drho}{\dot{\rho}}
\newcommand{\bveps}{\bolds\varepsilon}
\newcommand{\tbveps}{{\tilde{\bveps}}}
\newcommand{\tzeta}{\tilde{\zeta}}
\newcommand{\cor}{\operatorname{cor}}

\newtheorem{theorem}{Theorem}
\newtheorem{proposition}{Proposition}
\newproclaim{example}{Example}[section]

\newproclaim{Condition}{Condition}

\makeatother

\begin{document}
\begin{frontmatter}

\title{The sparse Laplacian shrinkage estimator for high-dimensional
regression}
\runtitle{Sparse Laplacian shrinkage estimator}

\begin{aug}
\author[A]{\fnms{Jian} \snm{Huang}\corref{}\thanksref{t1}\ead[label=e1]{jian-huang@uiowa.edu}},
\author[B]{\fnms{Shuangge} \snm{Ma}\thanksref{t2}\ead[label=e2]{shuangge.ma@yale.edu}},
\author[C]{\fnms{Hongzhe} \snm{Li}\thanksref{t3}\ead[label=e3]{hongzhe@upenn.edu}} and
\author[D]{\fnms{Cun-Hui} \snm{Zhang}\thanksref{t4}\ead[label=e4]{cunhui@stat.rutgers.edu}}
\runauthor{Huang, Ma, Li and Zhang}
\affiliation{University of Iowa, Yale University, University of Pennsylvania and~Rutgers University}
\address[A]{J. Huang\\
Department of Statistics\\
\quad and Actuarial Science, 241 SH \\
University of Iowa \\
Iowa City, Iowa 52242\\
USA\\
\printead{e1}} %adresu isvedimo komanda gale!
\address[B]{S. Ma \\
Division of Biostatistics \\
School of Public Health \\
Yale University \\
New Haven, Connecticut 06520\\
USA\\
\printead{e2}}
\address[C]{H. Li \\
Department of Biostatistics\\
\quad and Epidemiology \\
University of Pennsylvania School\\
\quad of Medicine \\
Philadelphia, Pennsylvania 19104\\
USA\\
\printead{e3}}
\address[D]{C.-H. Zhang \\
Department of Statistics\\
\quad and Biostatistics \\
Rutgers University\\
Piscataway, New Jersey 08854\hspace*{7.13pt}\\
USA\\
\printead{e4}}
\end{aug}

\thankstext{t1}{Supported in part by NIH Grants R01CA120988,
R01CA142774 and NSF Grant DMS-08-05670.}

\thankstext{t2}{Supported in part by NIH Grants R01CA120988,
R01CA142774, R03LM009754 and R03LM009828.}

\thankstext{t3}{Supported in part by NIH Grants
R01ES009911 and R01CA127334.}

\thankstext{t4}{Supported in part by NSF Grants DMS-06-04571,
DMS-08-04626 and NSA Grant MDS 904-02-1-0063.}

% HISTORY:
\received{\smonth{6} \syear{2010}}
\revised{\smonth{5} \syear{2011}}

% ABSTRACT
%
\begin{abstract}
We propose a new penalized method for variable selection and estimation
that explicitly incorporates the correlation patterns among predictors.
This method is based on a combination of the minimax concave penalty
and Laplacian quadratic associated with a graph as the penalty
function. We call it the sparse Laplacian shrinkage (SLS) method. The
SLS uses the minimax concave penalty for encouraging sparsity and
Laplacian quadratic penalty for promoting smoothness among coefficients
associated with the correlated predictors. The SLS has a generalized
grouping property with respect to the graph represented by the
Laplacian quadratic. We show that the SLS possesses an oracle property
in the sense that it is selection consistent and equal to the oracle
Laplacian shrinkage estimator with high probability. This result holds
in sparse, high-dimensional settings with $p \gg n$ under reasonable
conditions. We derive a coordinate descent algorithm for computing the
SLS estimates. Simulation studies are conducted to evaluate the
performance of the SLS method and a real data example is used to
illustrate its application.
\end{abstract}

% KEYWORDS
%

\end{frontmatter}

%s1 ###
\section{Introduction}
There has been much work on penalized methods for variable selection
and estimation in high-dimensional regression models. Several important
methods have been proposed. Examples include estimators based on the
bridge penalty [\citet{FrankF93}], the $\ell_1$ penalty or the
least absolute shrinkage and selection operator
[LASSO, \citet{Tibshirani96}, \citet{CDS}], the smoothly clipped absolute
deviation (SCAD) penalty [\citet{Fan97}, \citet{FanL}] and the minimum
concave penalty [MCP, \citet{Zhang10}]. These methods are able to do
estimation and automatic variable selection simultaneously and provide
a~computationally feasible way for variable selection in
high-dimensional settings. Much progress has been made in understanding
the theoretical properties of these methods. Efficient algorithms have
also been developed for implementing these methods.

A common feature of the methods mentioned above is the independence
between the penalty and the correlation among predictors. This can lead
to unsatisfactory selection results, especially in $p \gg n$ settings.
For example, as pointed out by \citet{ZouH}, the LASSO tends
to only select one variable among a group of highly
correlated variables; and its prediction performance may not be as good
as the ridge regression if there exists high correlation among
predictors. To overcome these limitations, \citet{ZouH}
proposed the elastic net (Enet) method, which uses a combination of the
$\ell_1$ and $\ell_2$ penalties. Selection properties of the Enet and
adaptive Enet have also been studied by Jia and Yu (\citeyear{JiaYu}) and
\citet{ZouZ}. \citet{BondelR} proposed the OSCAR (octagonal
shrinkage and clustering algorithm for regression) approach, which uses
a~combination of the $\ell_1$ norm and a pairwise $\ell_{\infty}$
norm for the coefficients. Huang et al. (\citeyear{HBMZa}) proposed the Mnet
method, which uses a combination of the MCP and $\ell_2$ penalties.
The Mnet estimator is equal to the oracle ridge estimator with high
probability under certain conditions. These methods are effective in
dealing with certain types of collinearity among predictors and has the
useful grouping property of selecting and dropping highly correlated
predictors together. Still, these combination penalties do not use any
specific information on the correlation pattern among the predictors.

\citet{LLa} proposed a network-constrained regularization
procedure for variable selection and estimation in linear regression
models, where the predictors are genomic data measured on genetic
networks. \citet{LLb} considered the general problem of regression
analysis when predictors are measured on an undirected graph, which is
assumed to be known a~priori. They called their method a
graph-constrained estimation procedure or GRACE. The GRACE penalty is
a combination of the $\ell_1$ penalty and a~penalty that is the
Laplacian quadratic associated with the graph. Because the GRACE uses
the $\ell_1$ penalty for selection and sparsity, it has the same
drawbacks as the Enet
discussed above. In addition, the full knowledge of the graphical
structure for the predictors is usually not available, especially in
high-dimensional problems. \citet{DayeJeng} proposed the weighted
fusion method, which also uses a combination of the $\ell_1$ penalty
and a quadratic form that can incorporate information among correlated
variables for estimation and variable selection. \citet{TutzU}
studied a form of correlation based penalty, which can be
considered a special case of the general quadratic penalty. But this
approach does not do variable selection. The authors proposed a
blockwise boosting procedure in combination with the correlation based
penalty for variable selection. \citet{HebiriGeer} studied
the theoretical properties of the smoothed-Lasso and other $\ell
_1+\ell_2$-penalized methods in $p \gg n$ models. \citet{PXS}
studied a grouped penalty based on the $L_{\gamma}$-norm for
$\gamma> 1$ that smoothes the regression coefficients over a network.
In particular, when $\gamma=2$ and after appropriate rescaling of the
regression coefficients, this group $L_{\gamma}$ penalty simplifies to
the group Lasso [Yuan and Lin (\citeyear{YuanL})] with the nodes in the network as
groups. This method is capable of group selection, but it does not do
individual variable selection. Also, because the group $L_{\gamma}$
penalty is convex for $\gamma> 1$, it does not lead to consistent
variable selection, even at the group level.

We propose a new penalized method for variable selection and estimation
in sparse, high-dimensional settings that takes into account certain
correlation patterns among predictors. We consider a combination of the
MCP and Laplacian quadratic as the penalty function. We call the
proposed approach the sparse Laplacian shrinkage (SLS) method.
The SLS uses the MCP to promote sparsity and Laplacian quadratic
penalty to encourage smoothness among coefficients associated with the
correlated predictors. An important advantage of the MCP over the $\ell
_1$ penalty is that it leads to estimators that are nearly unbiased and
achieve selection consistency under weaker conditions [\citet{Zhang10}].

The contributions of this paper are as follows.

\begin{itemize}
\item First, unlike the existing methods that use an $\ell_1$
penalty for selection and a ridge penalty or a general $\ell_2$
penalty for dealing with correlated predictors, we use
the MCP to achieve nearly unbiased selection and proposed a concrete
class of quadratics, the Laplacians, for incorporating correlation
patterns among predictors in a local fashion. In particular, we suggest
to employ the approaches for network analysis for specifying the
Laplacians. This provides an implementable strategy for incorporating
correlation structures in high-dimensional data analysis.

\item Second, we prove that the SLS estimator is sign consistent
and equal to the oracle Laplacian shrinkage estimator under reasonable
conditions. This result holds for a large class of Laplacian
quadratics. An important aspect of this result is that it allows the
number of predictors to be larger than the sample size. In contrast,
the works of \citet{DayeJeng} and \citet{TutzU} do not
contain such results in $p \gg n$ models. The selection consistency
result of \citet{HebiriGeer} requires certain strong
assumptions on the magnitude of the smallest regression coefficient
(their Assumption C) and on the correlation between important and
unimportant predictors (their Assumption D), in addition to a~variant
of the restricted eigenvalue condition (their Assumption B). In
comparison, our assumption involving the magnitude of the regression
coefficients is weaker and we use a sparse Riese condition instead of
imposing restriction on the correlations among predictors. In addition,
our selection results are stronger in that the SLS estimator is not
only sign consistent, but also equal to the oracle Laplacian shrinkage
estimator with high probability. In general, similar results are not
available with the use of the $\ell_1$ penalty.

\item Third, we show that the SLS method is potentially capable of
incorporating correlation structure in the analysis without incurring
extra bias. The Enet and the more general $\ell_1+\ell_2$ methods in
general introduces extra bias due to the quadratic penalty, in addition
to the bias resulting from the $\ell_1$ penalty. To the best of our
knowledge, this point has not been discussed in the existing
literature. We also demonstrate that
the SLS has certain local smoothing property with respect to the
graphical structure of the predictors.

\item Fourth, unlike in the GRACE method, the SLS does not assume
that the graphical structure for the predictors is known a priori. The
SLS uses the existing data to construct the graph Laplacian or to
augment partial knowledge of the graph structure.

\item Fifth, our simulation studies demonstrate that the SLS method
outperforms the $\ell_1$ penalty plus a quadratic penalty approach as
studied in \citet{DayeJeng} and \citet{HebiriGeer}. In
our simulation examples, the SLS in general has smaller empirical false
discovery rates with comparable false negative rates. It also has
smaller prediction errors.

\end{itemize}

This paper is organized as follows. In Section \ref{sec2}, we define the SLS
estimator. In Section \ref{sec3} we discuss ways to construct graph Laplacian,
or equivalently, its corresponding adjacency matrix. In Section \ref{sec4}, we
study the selection properties of the SLS estimators. In Section \ref{sec5}, we
investigate the properties of Laplacian shrinkage. In Section \ref{sec6}, we
describe a coordinate descent algorithm for computing the SLS
estimators, present simulation results and an application of the SLS
method to a microarray gene expression dataset. Discussions of the
proposed method and results are given in Section~\ref{sec7}. Proofs for the
oracle properties of the SLS and other technical details are provided
in the \hyperref[app]{Appendix}.

%s2 ###
\section{The sparse Laplacian shrinkage estimator}\label{sec2}

Consider the linear regression model
%
%e2.1 ###
\begin{equation}
\label{RegA} \by=\sum_{j=1}^p \bx_j \beta_j + \bveps
\end{equation}
with $n$ observations and $p$ potential predictors, where $\by=(y_1,
\ldots, y_n)'$ is the vector of $n$ response variables, $\bx
_{j}=(x_{1j}, \ldots, x_{nj})'$ is the $j$th predictor, $\beta_j$ is
the $j$th regression coefficient and $\bveps=(\veps_1, \ldots, \veps
_n)'$ is the vector of random errors.
Let $X=(\bx_1, \ldots, \bx_p)$ be the $n\times p$ design matrix.
Throughout, we assume that the response and predictors are centered and
the predictors are standardized
so that $\sum_{i=1}^n x_{ij}^2 = n,  j=1, \ldots,p$. For $\lam
=(\lam_1, \lam_2)$ with
$\lam_1 \ge0 $ and $\lam_2 \ge0$, we propose the penalized least
squares criterion
%
%e2.2 ###
\begin{eqnarray}
\label{SLSdefA}
M(\bb; \lam, \gam)
&=&\frac{1}{2n}\|\by-X\bb\|^2 + \sum_{j=1}^p\rho(|b_j|;\lam_1, \gam)
\nonumber\\[-8pt]\\[-8pt]
&&{} + \frac{1}{2} \lam_2
\sum_{1\le j < k \le p} \vert a_{jk}\vert(b_j-s_{jk}b_k)^2,\nonumber
\end{eqnarray}
where \mbox{$\|\cdot\|$} denotes the $\ell_2$ norm, $\rho$ is the MCP with
penalty parameter $\lam_1$ and regularization parameter $\gam$,
$|a_{jk}|$ measures the strength of the connection
between $\bx_j$ and $\bx_k$, and $s_{jk}=\sgn(a_{jk})$ is the sign
of $a_{jk}$, with $\sgn(t)=-1,0$ or $ 1$, respectively, for $t < 0,
=0$ or $> 0$. The two penalty terms in~(\ref{SLSdefA}) play different roles.
The first term promotes sparsity in the estimated model.
The second term encourages smoothness of the estimated coefficients of
the connected predictors. We can associate the quadratic form in this
term with the Laplacian for a suitably defined undirected weighted
graph for the predictors. See the description below.
For any given $(\lam, \gam)$, the SLS estimator~is
%
%e2.3 ###
\begin{equation}\label{SLS}
\hbbeta(\lam,\gam) = \argmin_{\bb} M(\bb; \lam,\gam).
\end{equation}

The SLS uses the MCP, defined as
%
%e2.4 ###
\begin{equation}
\label{MCPa}
\rho(t;\lam_1,\gam) = \lam_1 \int_0^{|t|} \bigl(1-x/(\gamma\lam
_1)\bigr)_{+} \,dx,
% t \ge0,
\end{equation}
where for any $a \in R$, $a_{+}$ is the nonnegative part of $a$, that
is, $a_+= a 1_{\{a \ge0\}}$. The MCP can be easily understood by
considering its derivative,
%
%e2.5 ###
\begin{equation}
\label{dMCPa}
\drho(t;\lam_1,\gam) = \lam_1 \bigl(1-|t|/(\gamma\lam_1)\bigr)_{+} \sgn(t).
\end{equation}
We observe that the MCP begins by applying the same level of
penalization as the $\ell_1$ penalty,
but continuously reduces that level to 0 for $|t| > \gamma\lambda$.
The regularization parameter $\gam$ controls the degree of concavity.
Larger values of $\gam$
make $\rho$ less concave. By sliding the value of $\gam$ from 1 to
$\infty$, the MCP provides a continuum of penalties with the
hard-threshold penalty as $\gam\rightarrow1+$ and the convex $\ell
_1$ penalty at $\gam= \infty$. Detailed discussion of MCP can be
found in \citet{Zhang10}.

The SLS also allows the use of different penalties than the MCP for
$\rho$, including the SCAD
[\citet{Fan97}, \citet{FanL}] and other quadratic splines. Because the
MCP minimizes the maximum concavity measure and has the simplest form\vadjust{\goodbreak}
among nearly unbiased penalties in this family, we choose it as the
default penalty for the SLS. Further discussion of the MCP and its
comparison with the LASSO and SCAD can be found in \citet{Zhang10} and
\citet{MFH}.

We express the nonnegative quadratic form in the second penalty term in~(\ref{SLSdefA})
using a positive semi-definite matrix $L$, which satisfies
\[
\bb'L\bb= \sum_{1 \le j < k \le p} \vert a_{jk}\vert
(b_j-s_{jk}b_k)^{2}\qquad \forall\bb\in\real^p.
\]
For simplicity, we confine our discussion to the symmetric case where
$a_{kj}=a_{jk}, 1 \le j <k \le p$. Since the diagonal elements $a_{jj}$
do not appear in the quadratic form, we can define them any way we like
for convenience. Let $A=(a_{jk}, 1\le j, k \le p)$ and $D=\diag(d_1,
\ldots, d_p)$, where $d_j=\sum_{k=1}^p \vert a_{jk}\vert$. We have
\(
\sum_{1 \le j < k \le p} |a_{jk}|(b_j-s_{jk}b_k)^2=\bb'(D-A)\bb.
\)
Therefore,\vspace*{1pt} $L=D-A$. This matrix is associated with a labeled
weighted graph $\cG=(V, \cE)$ with vertex set $V=\{1, \ldots, p\}$
and edge set $\cE=\{(j, k)\dvtx (j, k) \in V \times V\}$. Here the
$|a_{jk}|$ is the weight of edge $(j, k)$ and $d_j$ is the degree of
vertex $j$. The $d_j$ is also called the connectivity of vertex $j$.
The matrix $L$ is called the Laplacian of $\cG$ and $A$ its signed
adjacency matrix [\citet{Chunga}]. The edge $(j, k)$ is labeled with the
``$+$'' or ``$-$'' sign, but its weight $|a_{jk}|$ is always
nonnegative. We use a labeled graph to accommodate the
case where two predictors can have a nonzero adjacency coefficient but
are negatively correlated. Note that the usual adjacency matrix can be
considered a special case of signed adjacency matrix when all $a_{jk}
\ge0$. For simplicity, we will use the term adjacency matrix below.

We usually require that the adjacency matrix to be sparse in the sense
that many of its entries are zero or nearly zero. With a sparse
adjacency matrix, the main characteristic of the shrinkage induced by
the Laplacian penalty is that it occurs locally for the coefficients
associated with the predictors connected in the graph.
Intuitively, this can be seen by writing
\[
\lam_2 \sum_{1 \le j < k \le p} |a_{jk}|(b_j-s_{jk}b_k)^2 = \frac
{1}{2}\lam_2 \sum_{(j, k)\dvtx a_{jk}\neq0}|a_{jk}|(b_j-s_{jk}b_k)^2.
\]
Thus for $\lam_2 > 0$, the Laplacian penalty shrinks $b_j-s_{jk}b_k$
toward zero for $a_{jk}\neq0$. This can also be considered as a type
of local smoothing on the graph~$\cG$ associated with the adjacency
matrix $A$. In comparison, the shrinkage induced by the ridge penalty
used in the Enet is global in that
it shrinks all the coefficients toward zero, regardless of the
correlation structure among the predictors. We will discuss the
Laplacian shrinkage in more detail in Section \ref{sec5}.

Using the matrix notation, the SLS criterion (\ref{SLSdefA}) can be
written as
%
%e2.6 ###
\begin{equation}
\label{SLSdefB}\quad
M(\bb; \lam, \gam)
=\frac{1}{2n}\|\by-X\bb\|^2 + \sum_{j=1}^p\rho(|b_j|;\lam_1, \gam)
+ \frac{1}{2} \lam_2 \bb'(D-A)\bb.
\end{equation}
Here the Laplacian is not normalized, meaning that the weight $d_j$ is
not standardized to 1.
In problems where predictors should be treated without preference with
respect to connectivity, we can first normalized the Laplacian
$L^*=I_p-A^*$ with $A^*=D^{-1/2}AD^{-1/2}$ and use the criterion
\[
M^*(\bb;\lam,\gam) = \frac{1}{2n}\|\by-X\bb\|^2
+\sum_{j=1}^p\rho(|b_j|;\lam_1, \gam) + \frac{1}{2}\lam_2 \bb
^{\prime}(\bI_p-A^*)\bb.
\]
Technically, a normalized Laplacian $L^*$ can be considered a special
case of a~general $L$. We only consider the SLS estimator based on the
criterion~(\ref{SLSdefB}) when studying its properties.
In network analysis of gene expression data, genes with large
connectivity also tend to have important biological functions
[\citet{ZhangH05}]. Therefore, it is prudent to provide more
protection for such genes in the selection process.

%s3 ###
\section{Construction of adjacency matrix}\label{sec3}
In this section, we describe several simple forms of adjacency measures
proposed by \citet{ZhangH05}, which have have been successfully
used in network analysis of gene expression data. The adjacency measure
is often defined based on the notion of dissimilarity or similarity.

\begin{longlist}[(iii)]
\item[(i)]
A basic and widely used dissimilarity measure is the Euclidean
distance. Based on this distance, we can define adjacency coefficient as
\(
a_{jk}=\phi(\|\bx_j-\bx_k\|/\sqrt{n}),
\)
where $\phi\dvtx [0,\infty) \mapsto[0, \infty)$.
A simple adjacency function is the threshold function $\phi(x)=1\{x
\le2r\}$. Then
%
%e3.1 ###
\begin{equation}
\label{DAA}
a_{jk}= \cases{
1, &\quad if $\|\bx_j-\bx_k\|/\sqrt{n} \le2r$, \cr
0, &\quad if $\|\bx_j-\bx_k\|/\sqrt{n} >2r$.}
\end{equation}

It is convenient to express $a_{jk}$ in terms of the Pearson's
correlation coeffi\-cient $r_{jk}$ between $\bx_j$ and $\bx_k$, where
$r_{jk}=\bx_j'\bx_k/(\|\bx_j\| \|\bx_k\|)$. For predictors that
are standardized with $\|\bx_j\|^2=n, 1\le j\le p$, we have
\(\|\bx_j-\bx_k\|^2/n=\allowbreak 2-2 r_{jk}.\)
Thus in terms of correlation coefficients, we can write
\mbox{$a_{jk}=1\{r_{jk} > r\}$}. We determine the value of $r$ based on the Fisher
transformation $z_{jk} = 0.5\log((1+r_{jk})/(1-r_{jk}))$. If the
correlation between~$\bx_j$ and~$\bx_k$ is zero, $\sqrt{n-3} z_{jk}$
is approximately distributed as $N(0, 1)$. We can use this to determine
a threshold $c$ for
$\sqrt{n-3} z_{jk}$. The corresponding threshold for $r_{jk}$ is
$r=(\exp(2c/ \sqrt{n-3})-1)/(\exp(2c/ \sqrt{n-3})+1)$.

We note that here we use the Fisher transformation to change the scale
of the correlation coefficients from $[-1, 1]$ to the normal scale for
determining the threshold value $r$, so that the adjacency matrix is
relatively sparse. We are not trying to test the significance of
correlation coefficients.

\item[(ii)]
The adjacency coefficient in (\ref{DAA}) is defined based on a
dissimilarity measure. Adjacency coefficient can also be defined based
on similarity measures.
An often used similarity measure is Pearson's correlation
coefficient~$r_{jk}$.\vadjust{\goodbreak}
Other correlation measures such as Spearman's correlation can also be used. Let
\[
s_{jk}=\sgn(r_{jk}) \quad\mbox{and}\quad a_{jk}=s_{jk} 1\{|r_{jk}| > r\}.
\]
Here $r$ can be determined using the Fisher transformation as above.

\item[(iii)]
With the power adjacency function considered in \citet{ZhangH05},
\[
a_{jk}=\max(0, r_{jk})^{\alpha}\quad\mbox{and}\quad s_{jk}=1.
\]
Here $\alpha>0$ and can be determined by, for example, the scale-free
topology criterion.

\item[(iv)]
A variation of the above power adjacency function is
\[
a_{jk}=|r_{jk}|^{\alpha}\quad\mbox{and}\quad s_{jk}=\sgn(r_{jk}).
\]
\end{longlist}

For the adjacency matrices given above, (i) and (ii) use dichotomized
measures, whereas (iii) and (iv) use continuous measures. Under (i) and~(iii),
two covariates are either positively or not
connected/correlated. In contrast, under (ii) and (iv), two covariates
are allowed to be negatively connected/correlated.

There are many other ways for constructing an adjacency matrix.
For example, a popular adjacency measure in cluster analysis is
$a_{jk} = \exp(-\|\bx_j-\bx_k\|^2/n\tau^2)$ for $\tau> 0$. The
resulting adjacency matrix $A=[a_{jk}]$ is the Gram matrix associated
with the Gaussian kernel. For discrete covariates,
the Pearson correlation coefficient can still be used as a measure of
correlation or association between two discrete predictors or between a
discrete predictor and a continuous one. For example, for single
nucleotide polymorphism data, Pearson's correlation coefficient is
often used as a measure of linkage disequilibrium (i.e., association)
between two markers. Other measures, such as odds ratio or measure of
association based on contingency table can also be used for $r_{jk}$.

We note that how to construct the adjacency matrix is problem specific.
Different applications may require different adjacency matrices.
Since construction of adjacency matrix is not the focus
of the present paper, we will only consider the use of the four
adjacency matrices described above in our numerical studies in Section \ref{sec6}.

%s4 ###
\section{Oracle properties}
\label{sec4}
In this section, we study the theoretical properties of the SLS
estimator. Let the true value of the regression coefficient be
$\bbeta^o=(\beta_1^o, \ldots, \beta_p^o)'$. Denote
$\cO=\{j\dvtx \beta_j^o \neq0\}$, which is the set of indices of nonzero
coefficients. Let $d^o=|\cO|$ be the cardinality of $\cO$.
Define
%
%e4.1 ###
\begin{equation}
\label{ORA}
\hbbeta{}^o(\lam_2)=\argmin_{\bb}\biggl\{ \frac{1}{2n}\|\by-X\bb\|^2 +
\frac{1}{2}\lam_2\bb'L \bb, b_j =0, j \notin\cO\biggr\}.
\end{equation}
This is the oracle Laplacian shrinkage estimator on the set $\cO$.
Theorems \ref{ThmA} and \ref{ThmB} below provide sufficient\vadjust{\goodbreak}
conditions under which $\rP(\sgn(\hbbeta) \neq\sgn(\bbeta^o)$
or \mbox{$\hbbeta\neq\hbbeta{}^o) \rightarrow0$}. Thus, under those
conditions, the SLS estimator is sign consistent and equal to $\hbbeta
{}^o$ with high probability.

We need the following notation in stating our results.
Let $\Sigma= n^{-1} X'X$. For any $A\cup B \subseteq\{1, \ldots, p\}
$, vectors $\bv$, the design matrix $X$ and $V=(v_{ij})_{p\times p}$, define
\begin{eqnarray*}
\bv_B&=&(v_j, j\in B)',\qquad X_B = (\bx_j, j \in B),\\
V_{A,B}&=&(v_{ij}, i\in
A, j\in B)_{|A|\times|B|},\qquad
V_B=V_{B,B}.
\end{eqnarray*}
For example, $\Sigma_B = X_B'X_B/n$ and $\Sigma_{\cO}(\lam
_2)=\Sigma_{\cO}+\lam_2L_{\cO}$.
Let $|B|$ denote the cardinality of $B$.
Let $c_{\min}(\lam_2)$ be the smallest eigenvalue of $\Sigma+\lam_2 L$.
We use the following constants to bound the bias of the Laplacian:
%
%e4.2 ###
\begin{eqnarray}\label{C1C2}
C_1&=&\|\Sigma_{\mathcal{O}}^{-1}(\lam_2)L_\mathcal{O}\bbeta
^o_\mathcal{O}\|_\infty,\nonumber\\[-8pt]\\[-8pt]
C_2&=&\|\{\Sigma_{\cO^c,\mathcal{O}}(\lam_2)
\Sigma_\mathcal{O}^{-1}(\lam_2)L_\mathcal{O}-L_{\cO^c,\mathcal
{O}}\}\bbeta^o_\mathcal{O}\|_\infty.\nonumber
\end{eqnarray}

We make the following sub-Gaussian assumption on the error terms in~(\ref{RegA}).
\renewcommand{\theCondition}{(\Alph{Condition})}
\begin{Condition}\label{condA}
For a certain constant $\eps\in(0,1/3)$,
\[
\sup_{\|\bu\|=1}P\{\bu'\bveps> \sigma t\}\le e^{-t^2/2},\qquad
0<t\le\sqrt{2\log(p/\eps)}.
\]
\end{Condition}

%s4.1 ###
\subsection{Convex penalized loss}\label{sec41}
We first consider the case where $\Sigma(\lam_2)=\Sigma+\lam_2L$ is
positive definite. Since (\ref{ORA}) is the minimizer of the Laplacian
restricted to the support $\mathcal{O}$,
it can be explicitly written as
%
%e4.3 ###
\begin{equation}
\label{KKTa}
\hbbeta{}^o_\mathcal{O} = (\Sigma_\mathcal{O}+\lam_2L_\mathcal
{O})^{-1}X_\mathcal{O}'\by/n,\qquad \hbbeta{}^o_{\mathcal{O}^c}=0,
\end{equation}
provided that $\Sigma_\mathcal{O}(\lam_2)$ is invertible. Its
expectation $\bbeta^*=E\hbbeta{}^o$,
considered as a target of the SLS estimator, must satisfy
%
%e4.4 ###
\begin{equation}\label{target}
\bbeta^*_\mathcal{O} = (\Sigma_\mathcal{O}+\lam_2L_\mathcal
{O})^{-1}\Sigma_\mathcal{O}\bbeta^o,\qquad
\bbeta^*_{\mathcal{O}^c}=0.
\end{equation}
\begin{Condition}\label{condB}
(i) $c_{\min}(\lam_2)>1/\gamma$ with $\rho
(t;\lam_1,\gamma)$
in (\ref{SLSdefA}).

\mbox{\hphantom{i}}(ii) The penalty levels satisfy
\[
\lam_1\ge\lam_2C_2 + \sigma\sqrt{2\log\bigl((p-d^o)/\eps\bigr)}\max_{j\le
p}\|\bx_j\|/n
\]
with $C_2$ in (\ref{C1C2}).

(iii) With $\{v_j,j\in\mathcal{O}\}$ being the diagonal elements of
$\Sigma_\mathcal{O}^{-1}(\lam_2)\Sigma_\mathcal{O}\{\Sigma
_\mathcal{O}^{-1}(\lam_2)\}$,
\[
\min_{j\in\mathcal{O}}\{|\beta^*_j|(n/v_j)^{1/2}\}
\ge\sigma\sqrt{2\log(d^o/\eps)}.
\]
\end{Condition}

Define $\beta_*= \min\{|\beta_j^o|, j\in\cO\}$. If $\cO$ is an
empty set, that is, when all the regression coefficients are zero,
we set $\beta_*=\infty$.\vadjust{\goodbreak}
\begin{theorem}
\label{ThmA}
Suppose Conditions \ref{condA} and \ref{condB} hold. Then
%
%e4.5 ###
\begin{equation}\label{ThmA-1}
\rP(\{j\dvtx \hbeta_j\neq0\}\neq\cO\mbox{ or } \hbbeta\neq
\hbbeta{}^o) \le3\eps.
\end{equation}
If $\beta_*
\ge\lam_2C_1+\max_j\sqrt{(2v_j/n)\log(d^o/\eps)}$ instead of
Condition \ref{condB}\textup{(iii)}, then
%
%e4.6 ###
\begin{equation}\label{ThmA-2}
\rP\bigl(\sgn(\hbbeta) \neq\sgn(\bbeta^o) \mbox{ or } \hbbeta
\neq\hbbeta{}^o\bigr) \le3\eps.
\end{equation}
Here note that $p, d^o, \gamma$ and $c_{\min}(\lam_2)$ are all
allowed to depend on $n$.
\end{theorem}

The probability bound on the selection error in Theorem \ref{ThmA} is
nonasymptotic. If the conditions of Theorem \ref{ThmA} hold with $\eps
\to0$, then
(\ref{ThmA-1}) implies selection consistency of the SLS estimator
and (\ref{ThmA-2}) implies sign consistency.
The conditions are mild. Condition \ref{condA} concerns the tail probabilities
of the error distribution and is satisfied if the errors are normally
distributed. Condition \ref{condB}(i) ensures that the SLS criterion is
strictly convex so that the solution is unique.
The oracle estimator $\hbeta{}^o$ is biased due to the Laplacian
shrinkage. Condition \ref{condB}(ii) requires a penalty level $\lam_1$ to
prevent this bias and noise to cause false selection of variables in
$\cO^c$. Condition \ref{condB}(iii) requires that the nonzero coefficients
not be too small in order for the SLS estimator to be able to
distinguish nonzero from zero coefficients.

In Theorem \ref{ThmA}, we only require $c_{\min}(\lam_2) > 0$, or
equivalently, $\Sigma+\lam_2L$~to be positive definite. The matrix
$\Sigma$ can be singular. This can be seen~as follows.
The adjacency matrix partitions the graph into
disconnected cliques~$V_g$, $1 \le g \le J$, for some $J \ge1$.
Let node $j_g$ be a (representative) member of $V_g$. A~node $k$
belongs to the same clique $V_g$ iff (if and only if)\vspace*{1pt}
$a_{j_gk_1}a_{k_1k_2}\cdots a_{k_m k}\neq0$ through\vspace*{1pt} a certain chain
$j_g\to k_1\to k_2\to\cdots\to k_m\to k$. Define
${\bar\bx}_g=\break\sum_{k\in V_g}a_{j_gk_1}a_{k_1k_2}\cdots a_{k_m k}\bx
_k/|V_g|$,\vspace*{1pt}
where $|V_g|$ is the cardinality of~$V_g$.
The matrix $\Sigma+\lam_2L$ is positive definite iff $\bb'\Sigma\bb
=\bb'L\bb=0$
implies $\bb=0$. Since $\bb'L\bb=0$ implies $\sum_{k\in V_g}b_k\bx
_k=b_{j_g}|V_g|{\bar\bx}_g$,
$\Sigma+\lam_2L$ is positive definite iff the vectors ${\bar\bx}_g$
are linearly independent. This does not require $n\ge p$. In other
words, Theorem \ref{ThmA} is applicable to $p > n$ problems as long as
the vectors ${\bar\bx}_g$ are linearly independent.

%s4.2 ###
\subsection{The nonconvex case}
When $\Sigma(\lam_2)=\Sigma+\lam_2 L$ is singular, Theorem \ref
{ThmA} is not applicable. In this case, further conditions are required
for the oracle property to hold. The key condition needed is the sparse
Reisz condition, or SRC [\citet{ZhangHuang08}], in (\ref{SRCa}) below.
It restricts the spectrum of diagonal subblocks of $\Sigma(\lam_2)$
up to a certain dimension.\vspace*{1pt}

Let $\tX=\tX(\lam_2)$ be a matrix satisfying $\tX'\tX/n=\Sigma
(\lam_2)=X'X/n+\lam_2L$ and
$\tby= \tby(\lam_2)$ be a vector satisfying $\tX'\tby= X'\by$. Define
%
%e4.7 ###
\begin{equation}
\label{tM}
\tM(\bb;\lam,\gamma) = \frac{1}{2n}\|\tby- \tX\bb\|^2+\sum
_{j=1}^p\rho(|b_j|;\lam_1,\gamma).
\end{equation}
Since $M(\bb;\lam,\gamma) - \tM(\bb;\lam,\gamma) = (\|\by\|^2-\|
\tby\|^2)/(2n)$,
the two penalized loss functions have the same set of local minimizers.
For the penalized loss (\ref{tM}) with the data $(\tX,\tby)$, let
%
%e4.8 ###
\begin{equation}
\label{MC+}
\hbbeta(\lam) = \bdelta(\tX(\lam_2),\tby(\lam_2),\lam
_1),
\end{equation}
where the map $\bdelta(X,\by,\lam_1)\in\real^p$ defines the MC$+$
estimator [\citet{Zhang10}] with data $(X,\by)$ and penalty level $\lam
_1$. It was\vspace*{1pt} shown in \citet{Zhang10} that $\bdelta(X,\by,\lam_1)$
depends on $(X,\by)$ only through $X'\by/n$ and $X'X/n$, so that
different choices of $\tX$ and $\tby$ are allowed. One way\vspace*{1pt} is to pick
$\tby=(\by',0)'$ and
$\tX=\diag(X,(n\lam_2L)^{1/2})$. Another way is to pick
$\tX'\tX/n=\Sigma(\lam_2)$ and $\tby= (\tX')^\dag X'\by$ of
smaller dimensions, where $(\tX')^\dag$ is the Moore--Penrose inverse
of $\tX'$.
\begin{Condition}\label{condC}
(i) For an integer $d^*$ and spectrum bounds
$0 < c_*(\lam_2) \le c^*(\lam_2) < \infty$,
%
%e4.9 ###
\begin{eqnarray}
\label{SRCa}
0 < c_*(\lam_2) \le\bu_B'\Sigma_B(\lam_2)\bu_B \le c^*(\lam_2)<
\infty\nonumber\\[-8pt]\\[-8pt]
&&\eqntext{\forall B \mbox{ with } |B\cup\cO| \le d^*, \|\bu_B\|=1,}
\end{eqnarray}
with $d^* \ge d^o(K_*+1)$, $\gamma\ge c_*^{-1}(\lam_2)\sqrt
{4+c_*(\lam_2)/c^*(\lam_2)}$
in (\ref{SLSdefA}), and $K_*=c^*(\lam_2)/c_*(\lam_2)-1/2$.

\mbox{\hphantom{i}}(ii) With $C_2=\|\{\Sigma_{B,\mathcal{O}}(\lam_2)
\Sigma_\mathcal{O}^{-1}(\lam_2)L_\mathcal{O}-L_{B,\mathcal{O}}\}
\bbeta^o_\mathcal{O}\|_\infty$,
\[
\max\bigl\{1,\sqrt{c_*(\lam_2)K_*/(K_*+1)}\bigr\}\lam_1\ge\lam_2C_2
+ \sigma\sqrt{2\log(p/\eps)}\max_{j\le p}\|\bx_j\|/n.
\]

(iii) With $\{v_j,j\in\mathcal{O}\}$ being the diagonal
elements of
$\Sigma_\mathcal{O}^{-1}(\lam_2)\Sigma_\mathcal{O}\{\Sigma
_\mathcal{O}^{-1}(\lam_2)\}$,
\[
\min_{j\in\mathcal{O}}\bigl\{|\beta^*_j|-\gamma\bigl(2\sqrt{c^*(\lam
_2)}\lam_1\bigr)\bigr\}(n/v_j)^{1/2}
\ge\sigma\sqrt{2\log(d^o/\eps)}.
\]
\end{Condition}
\begin{theorem}
\label{ThmB}
\textup{(i)} Suppose Conditions \ref{condA} and \ref{condC} hold. Let $\hbbeta(\lam)$ be as in
(\ref{MC+}). Then
%
%e4.10 ###
\begin{equation}\label{ThmB-1}%\begin{eqnarray*}
\rP(\{j\dvtx \hbeta_j\neq0\}\neq\cO\mbox{ or } \hbbeta\neq
\hbbeta{}^o) \le3\eps.
\end{equation}
If $\beta_*
\ge\lam_2C_1+\gamma(2\sqrt{c^*(\lam_2)}\lam_1)+
\max_j\sqrt{(2v_j/n)\log(d^o/\eps)}$ instead of
Condition~\ref{condC}\textup{(iii)}, then
%
%e4.11 ###
\begin{equation}\label{ThmB-2}%\begin{eqnarray*}
\rP\bigl(\sgn(\hbbeta) \neq\sgn(\bbeta^o) \mbox{ or } \hbbeta
\neq\hbbeta{}^o\bigr) \le3\eps.
\end{equation}
Here note that $p$, $\gamma$, $d^o$, $d^*$, $K_*$, $\eps$, $c_*(\lam
_2)$ and $c^*(\lam_2)$ are all allowed to depend on~$n$,
including the case $c_*(\lam_2) \rightarrow0$ as long as the
conditions hold as stated.

\textup{(ii)} The statements in \textup{(i)} also hold for all local minimizers $\hbeta$
of (\ref{SLSdefB}) or~(\ref{tM})
satisfying $\#\{j\notin\cO\dvtx \hbeta_j\neq0\}+d^o \le d^*$.
\end{theorem}

If the conditions of Theorem \ref{ThmB} hold with $\eps\to0$, then
(\ref{ThmB-1}) implies selection consistency of the SLS estimator
and (\ref{ThmB-2}) implies sign consistency.\vadjust{\goodbreak}

Condition \ref{condC}, designed to handle the noncovexity of the penalized
loss, is a weaker version of Condition \ref{condB} in the sense of allowing
singular $\Sigma(\lam_2)$. The SRC (\ref{SRCa}), depending on $X$ or
$\tX$ only through the regularized Gram matrix $\tX'\tX/n=\Sigma
(\lam_2)=\Sigma+\lam_2L$,
ensures that the model is identifiable in a lower $d^*$-dimensional
space. When $p > n$, the smallest singular value of $X$ is always zero.
However, the requirement $c_*(\lam_2)> 0$ only concerns $d^*\times
d^*$ diagonal submatrices of $\Sigma(\lam_2)$, not the Gram matrix
$\Sigma$ of the design matrix $X$. We can have $p\gg n$ but still
require $d^*/d^o\ge K_*+1$ as in (\ref{SRCa}). Since $p, d^0, \gam$,
$d^*$, $K_*$, $c_*(\lam_2)$ and $c^*(\lam_2)$ can depend on $n$, we
allow the case $c_*(\lam_2) \rightarrow0$ as long as Conditions \ref{condA}
and \ref{condC} hold as stated. Thus, we allow $p \gg n$ but require that the
model is sparse, in the sense that the number of nonzero coefficients
$d^o$ is smaller than $d^*/(1+K_*)$. For example, if $c_*(\lam
_2)\asymp O(n^{-\alpha})$ for a small\vspace*{1pt} $\alpha> 0$ and $c^*(\lam
_2)\asymp O(1)$, then we require $\gamma\asymp O(n^{3\alpha/2})$ or
greater, $K^*\asymp O(n^{\alpha})$ and $d^*/d^o \asymp O(n^{\alpha})$
or greater. So all these quantities can depend on $n$, as long as the
other requirements are met in Condition~\ref{condC}.

By examining the Conditions \ref{condC}(ii) and \ref{condC}(iii),
for standardized predictors with $\|\bx_j\|=\sqrt{n}$, we can have
$\log(p/\eps)=o(n)$ or $p=\eps\exp(o(n))$ as long as Condition
\ref{condC}(ii) is satisfied. As in \citet{Zhang10}, under a somewhat stronger
version of Condition \ref{condC}, Theorem \ref{ThmB} can be extended to
quadratic spline concave penalties satisfying $\rho(t;\lam_1,\gamma
)=\lam_1^2\rho(t/\lam;\gamma)$ with a penalty function satisfying
$(\partial/\partial t)\rho(t;\gamma)=1$ at $t=0+$ and $0$ for
$t>\gamma$.

Also, comparing our results with the selection consistency results of
\citet{HebiriGeer} on the smoothed $\ell_1 + \ell
_2$-penalized methods, our conditions tend to be weaker. Notably,
\citet{HebiriGeer} require an condition on the Gram matrix
which assumes that the correlations between the truly relevant
variables and those which are not are small. No such assumption is
required for our selection consistency results. In addition, our
selection results are stronger in the sense that the SLS estimator is
not only sign consistent, but also equal to the oracle Laplacian
shrinkage estimator with high probability. In general, similar results
are not available with the use of the $\ell_1$ penalty for sparsity.

Theorem \ref{ThmB} shows that the SLS estimator automatically adapts
to the sparseness of the $p$-dimensional model and the denseness of a
\textit{true} submodel. From a sparse $p$-model, it correctly selects the
true underlying model $\cO$. This underlying model is a dense model in
the sense that all its coefficients are nonzero. In this dense model,
the SLS estimator behaves like the oracle Laplacian shrinkage estimator
in (\ref{ORA}). As in the convex penalized loss setting, here the
results do not require a correct specification of
a population correlation structure of the predictors.

%s4.3 ###
\subsection{Unbiased Laplacian and variance reduction}
There are two natural questions concerning the SLS. First, what are the
benefits from introducing the Laplacian penalty? Second,\vadjust{\goodbreak}
what kind of Laplacian $L$ constitutes a reasonable choice?
Since the SLS estimator is equal to the oracle Laplacian estimator with
high probability by Theorem \ref{ThmA} or \ref{ThmB}, these questions
can be answered by examining the oracle Laplacian shrinkage estimator
(\ref{ORA}), whose nonzero part is
\[
\hbbeta{}^o_{\cO}(\lam_2) = \Sigma_{\cO}^{-1}(\lam_2)X_{\cO}'\by/n.
\]
Without the Laplacian, that is, when $\lam_2=0$, it becomes the
least squares (LS) estimator
\[
\hbbeta{}^o_{\cO}(0) = \Sigma_{\cO}^{-1}X_{\cO}'\by/n.
\]
If some of the predictors in $\{\bx_j, j \in\cO\}$ are highly
correlated or $|\cO| \ge n$, the LS estimator $\hbbeta{}^o_{\cO}(0)$
is not stable or unique. In comparison, as discussed below Theorem~\ref
{ThmA}, $\Sigma_{\cO}(\lam_2) = \Sigma_{\cO}+\lam_2L_{\cO}$ can
be a full rank matrix under a~reasonable condition, even if the
predictors in $\{\bx_j, j \in\cO\}$ are highly correlated or $|\cO|
\ge n$.

For the second question, we examine the bias of $\hbbeta{}^o_{\cO}(\lam_2)$.
Since the bias of the target vector (\ref{target}) is $\bbeta{}^o_{\cO
}-\bbeta_{\cO}^*(\lam_2)=
\lam_2\Sigma_{\cO}^{-1}(\lam_2)L_{\cO}\bbeta_{\cO}^o$,
$\hbbeta{}^o_{\cO}(\lam_2)$ is unbiased iff $L_{\cO}\bbeta_{\cO}^o=0$.
Therefore, in terms of bias reduction, a Laplacian $L$ is most
appropriate if the condition $L_{\cO}\bbeta_{\cO}^o=0$ is satisfied.
We shall say that a Laplacian $L$ is unbiased if $L_{\cO}\bbeta_{\cO}^o=0$.
It follows from the discussion at the end of Section \ref{sec41} that
$L_{\cO}\bbeta_{\cO}^o=0$ if $\beta^o_k = \beta
^o_{j_g}a_{j_gk_1}a_{k_1k_2}\cdots a_{k_m k}$,
where $j_g$ is a representative member of the clique $V_g\cap\cO$ and
$\{k_1,\ldots,k_m,k\}\subseteq V_g\cap\cO$.

With an unbiased Laplacian, the mean square error of $\hbbeta{}^o_{\cO
}(\lam_2)$ is
\[
E\|\hbbeta{}^o_{\cO}(\lam_2)-\bbeta_{\cO}^o\|^2 = \frac{\sigma
^2}{n}\trace(\Sigma_{\cO}^{-1}(\lam_2)\Sigma_{\cO}\Sigma_{\cO
}^{-1}(\lam_2)).
\]
The mean square error of $\hbbeta_{\cO}(0)$ is
\[
E\|\hbbeta{}^o_{\cO}(0)-\bbeta_{\cO}^o\|^2 = \frac{\sigma
^2}{n}\trace(\Sigma_{\cO}^{-1}).
\]
We always have $E\|\hbbeta{}^o_{\cO}(\lam_2)-\bbeta_{\cO}^o\|^2
< E\|\hbbeta{}^o_{\cO}(0)-\bbeta_{\cO}^o\|^2$ for $\lam_2 > 0$.
Therefore, an unbiased Laplacian reduces variance without incurring any
bias on the estimator.

%s5 ###
\section{Laplacian shrinkage}\label{sec5}
The results in Section \ref{sec4} show that the SLS estimator is
equal to the oracle Laplacian shrinkage estimator with probability
tending to one under certain conditions. In addition, an unbiased
Laplacian reduces variance but does not increase bias.
Therefore, to study the shrinkage effect of the Laplacian penalty on
$\hbbeta$, we can consider the oracle estimator $\hbbeta{}^o_{\cO}$.
To simplify the notation and without causing confusion,
in this section,
we study some other basic properties of the Laplacian shrinkage and
compare it with the ridge shrinkage. The Laplacian shrinkage
estimator is defined as
%
%e5.1 ###
\begin{equation}
\qquad\tbbeta(\lam_2)=\argmin_{\bb}\biggl\{ G(\bb;\lam_2) \equiv\frac
{1}{2n} \|\by-X\bb\|^2
+ \frac{1}{2}\lam_2\bb'L\bb,  \bb\in\real^q\biggr\}.
\label{ORB}
\end{equation}

The following proposition shows that the Laplacian penalty shrinks
a~coefficient toward the center of all the coefficients connected to
it.\vspace*{3pt}

\begin{proposition}
\label{PropA}
Let $\tbr=\by-X\tbbeta$.

\textup{\mbox{}\hphantom{i}(i)}
\[
\lam_2 \max_{1\le j \le q} d_j |\tbeta_j-\ba_j'\tbbeta/d_j|
\le\|\tbr\| \le\|\by\|.
\]

\textup{(ii)}
\[
\lam_2|d_j\tbeta_j-\ba_j'\tbbeta- (d_k\tbeta_k-\ba_k'\tbbeta)|
\le\frac{1}{n} \|\bx_j-\bx_k\| \|\by\|.
\]
\end{proposition}

Note that $\ba_j'\tbbeta/d_j=\sum_{k=1}^q a_{jk}\tbeta_k/d_j
=\sum_{k=1}^q \sgn(a_{jk})|a_{jk}|\tbeta_k/d_j$ is a signed weighted
average of the $\tbeta_k$'s connected to $\tbeta_j$, since $d_j=\sum
_k |a_{jk}|$. Part (i) of Proposition \ref{PropA} provides an upper
bound on the difference between $\tbeta_j$ and the center of all the
coefficients connected to it. When $\|\tbr\|/(\lam_2d_j) \to0$, this
difference converges to zero. For standardized $d_j=1$,
part (ii) implies that the difference between the centered $\tbeta_j$
and $\tbeta_k$ converges to zero
when $\|\bx_j-\bx_k\| \|\by\|/(\lam_2 n) \rightarrow0$.

When there are certain local structures
in the adjacency matrix $A$, shrinkage occurs at the local level.
As an example, we consider the adjacency matrix based on partition of
the predictors into $2r$-balls defined in (\ref{DAA}). Correspondingly,
the index set $\{1, \ldots, q\}$ is divided into disjoint
neighborhoods/cliques $V_1, \ldots, V_J$. We consider the normalized
Laplacian $L=I_q - A$, where $I_q$ is a $q\times q $ identity matrix
and $A=\diag(A_1, \ldots, A_J)$ with
$A_g=v_g^{-1} \bone_g'\bone$. Here $v_g=|V_g|, 1 \le g \le J$.
Let $\bb_g=(b_j, j \in V_g)'$. We can write the objective function as
%
%e5.2 ###
\begin{equation}
\label{Na}
G(\bb;\lam_2)
= \frac{1}{2n}\|\by-X\bb\|^2
+ \frac{1}{2}\lam_2 \sum_{g=1}^J \bb_g'(\bI_g-v_g^{-1}\bone
_g'\bone_g)\bb_g.
\end{equation}
For the Laplacian shrinkage estimator based on this criterion, we have
the following grouping properties.\vspace*{3pt}

\begin{proposition}
\label{PropB}
\textup{(i)} For any $j, k \in V_g, 1 \le g \le J$,
\[
\lam_2\vert\tbeta_j-\tbeta_k\vert \le\frac{1}{n}\|\bx_j-\bx_k\|
\cdot\|\by\|, \qquad  j, k \in V_g.
\]

\textup{(ii)} Let $\bar{\beta}_g$ be the average of the estimates in $V_g$.
For any $j \in V_g$ and $k \in V_h$, $g\neq h$,
\[
\lam_2|\tbeta_j-\bar{\beta}_g-(\tbeta_k-\bar{\beta}_h)| \le
\frac{1}{n}\|\bx_j-\bx_k\|\cdot\|\by\|, \qquad  j \in V_g, k \in V_h.
\]
\end{proposition}

This proposition characterizes the smoothing effect
and grouping property of the Laplacian penalty in (\ref{Na}).
Consider the case $\|\by\|^2/n=O(1)$.
Part (i) implies that, for $j$ and $k$ in the same neighborhood and
$\lam_2 > 0$, the difference
$\tbeta_j-\tbeta_k \rightarrow0$ if $\|\bx_j-\bx_k\|/(\lam
_2n^{1/2}) \rightarrow0$. Part (ii) implies that, for
$j$ and $k$ in different neighborhoods and $\lam_2 > 0$, the
difference between the centered $\tbeta_j$ and $\tbeta_k$ converges
to zero if
$\|\bx_j-\bx_k\|/(\lam_2n^{1/2}) \rightarrow0$.

We now compare the Laplacian shrinkage and ridge shrinkage.
The discussion at the end of Section \ref{sec4} about the requirement for the
unbiasedness of Laplacian can be put in a wider context when a general
positive definite or semidefinite matrix $Q$ is used in the place of
$L$. This wider context includes the Laplacian shrinkage and ridge
shrinkage as special cases. Specifically, let
\[
\hbbeta_{Q}(\lam,\gam) =\argmin_{\bb} \frac{1}{2n}\|\by-X\bb\|
^2 + \sum_{j=1}^p\rho(|b_j|;\lam_1,\gam) + \frac{1}{2}\lam_2\bb
'Q\bb.
\]
For $Q=I_p$, $\hbbeta_Q$ becomes the Mnet estimator [Huang et al. (\citeyear{HBMZa})].
With some modifications on the conditions in Theorem \ref{ThmA} or
Theorem \ref{ThmB}, it can be shown that
$\hbbeta_{Q}$ is equal to the oracle estimator defined as
\[
\hbbeta{}^o_{Q}(\lam_2) =\argmin_{\bb} \biggl\{\frac{1}{2n}\|\by
-X\bb\|^2 + \frac{1}{2}\bb'Q\bb, b_j=0, j \notin\cO\biggr\}.
\]
Then in a way similar to the discussion in Section \ref{sec4}, $\hbbeta_{Q}$ is
nearly unbiased iff
$Q_{\cO}\bbeta_{\cO}^o=0$. Therefore, for $
\|\bbeta^o_{\cO}\|\neq0$, $Q_{\cO}$
must be a rank deficient matrix, which in turn implies that $Q$ must be
rank deficient. Note that any Laplacian $L$ is rank deficient. This
rank deficiency requirement excludes the ridge penalty with \mbox{$Q=I_p$}.
For the ridge penalty to yield an unbiased estimator, it must hold that
$\|\bbeta^o\|=0$ in the underlying model.

We now give a simple example that illustrates the basic characteristics
of Laplacian shrinkage and its differences from ridge
shrinkage.

\begin{example}
\label{AshrinkA}
\begin{rm}
Consider a linear regression model with two predictors satisfying $\|
\bx_j\|^2=n$, $j=1,2$.
The Laplacian shrinkage and ridge estimators are defined as
\[
(\hb_{L1}(\lam_2), \hb_{L2}(\lam_2)) = \argmin_{b_1, b_2} \frac
{1}{2n} \sum_{i=1}^n (y_i-x_{i1}b_1-x_{i2}b_2)^2 %+ \sum_{j=1}^2
+ \frac{1}{2}\lam_2 (b_1-b_2)^2
\]
and
\[
(\hb_{R1}(\lam_2), \hb_{R2}(\lam_2))= \argmin_{b1, b_2}\frac
{1}{2n} \sum_{i=1}^n (y_i-x_{i1}b_1-x_{i2}b_2)^2
+ \frac{1}{2}\lam_2 (b_1^2+b_2^2).\vadjust{\goodbreak}
\]
Denote $r_1=\cor(\bx_1, \by)$, $r_2=\cor(\bx_2, \by)$ and
$r_{12}=\cor(\bx_1, \bx_2)$.
The Laplacian shrinkage estimates are
\[
\hb_{L1}(\lam_2)= \frac{(1+\lam_2)r_1-(r_{12}-\lam_2)r_2}{(1+\lam
_2)^2-(r_{12}-\lam_2)^2},\qquad
\hb_{L2}(\lam_2)=
\frac{(1+\lam_2)r_2-(r_{12}-\lam_2)r_1}{(1+\lam_2)^2-(r_{12}-\lam_2)^2}.
\]
Let
\[
\hb_{\mathrm{ols}1}=\frac{r_1-r_{12}r_2}{1-r_{12}^2},\qquad
\hb_{\mathrm{ols}2}=\frac{r_2-r_{12}r_1}{1-r_{12}^2},\qquad
\hb_L(\infty) = \frac{r_1+r_2}{2(1+r_{12})},
\]
where $(\hb_{\mathrm{ols}1}, \hb_{\mathrm{ols}2})$ is the ordinary least squares (OLS)
estimator for the bivariate regression, $\hb_L(\infty)$ is the OLS
estimator that assumes the two coefficients are equal, that is, it
minimizes $ \sum_{i=1}^n(y_i-(x_{i1}+x_{i2})b)^2$.\vspace*{1pt}
Let $w_{L}=(2\lam_2)/(1-r_{12}+2\lam_2)$. After some simple algebra, we have
\[
\hb_{L1}(\lam_2)= (1-w_L) \hb_{\mathrm{ols}1} +
w_L \hb_L(\infty)
\]
and
\[
\hb_{L2}(\lam_2)= (1-w_L)\hb_{\mathrm{ols}2} + w_L \hb_L(\infty).
\]
Thus, for any fixed $\lam_2$, $\hb_L(\lam_2)$ is a weighted average
of $\hb_{\mathrm{ols}}$ and $\hb_L(\infty)$ with the weights depending on
$\lam_2$.
When $\lam_2 \rightarrow\infty$,
\(
\hb_{L1} \rightarrow\hb_L(\infty) \mbox{ and } \hb_{L2}
\rightarrow\hb_L(\infty).
\)
Therefore, the Laplacian penalty shrinks the OLS estimates toward
a common value, which is the OLS estimate assuming equal regression
coefficients.

Now consider the ridge regression estimator. We have
\[
\hb_{R1}(\lam_2)= \frac{(1+\lam_2)r_1-r_{12}r_2}{(1+\lam_2)^2-r_{12}^2}
\quad \mbox{and}\quad
\hb_{R2}(\lam_2)= \frac{(1+\lam_2)r_2-r_{12}r_1}{(1+\lam_2)^2-r_{12}^2}.
\]
The ridge estimator converges to zero as $\lam_2 \rightarrow\infty$.
For it to converge to a~nontrivial solution, we need to rescale it by a
factor of $1+\lam_2$.
Let $w_R=\lam/(1+\lam-r_{12}^2)$.
Let $\hb_{u1}=r_1$ and $\hb_{u2}=r_2$. Because $n^{-1}\sum
_{i=1}^nx_{i1}^2=1$ and $n^{-1}\sum_{i=1}^nx_{i2}^2=1$, $r_1$ and
$r_2$ are also the OLS estimators of univariate regressions of
$\by$ on $\bx_1$ and $\by$ on $\bx_2$, respectively.
We can write
\begin{eqnarray*}
(1+\lam_2)\hb_{R1}(\lam_2) &=& c_{\lam_2} (1-w_R)\hb_{\mathrm{ols}1}
+ c_{\lam} w_R \hb_{u1} ,
\\
(1+\lam_2)\hb_{R2}(\lam_2) &=& c_{\lam_2} (1-w_R) \hb_{\mathrm{ols}2}
+ c_{\lam} w_R \hb_{u2},
\end{eqnarray*}
where\vspace*{1pt} $c_{\lam_2} = \{(1+\lam_2)^2-(1+\lam)r_{12}^2\}/\{(1+\lam
_2)^2-r_{12}^2\}$.
Note that $c_{\lam_2} \approx1$. Thus, $(1+\lam_2)\hb_R$ is a
weighted average of
the OLS and the univariate regression estimators. The ridge penalty
shrinks the (rescaled) ridge estimates toward individual univariate
regression estimates.
\end{rm}
\end{example}

%s6 ###
\section{Simulation studies}\label{sec6}

We use a coordinate descent algorithm to compute the SLS estimate. This
algorithm optimizes a target function with respect to a single
parameter at a time and iteratively cycles through all parameters until\vadjust{\goodbreak}
convergence. This algorithm was originally proposed for criterions with
convex penalties such as LASSO [\citet{Fu}, \citet{GenkinLewisMadigan04},
\citet{FHHT}, Wu and Lange (\citeyear{WuL})]. It has been proposed to
calculate the MCP estimates [\citet{BH}]. Detailed steps of
this algorithm for computing the SLS estimates can be found in the
technical report accompanying this paper [Huang et al. (\citeyear{HBMZb})].

In simulation studies, we consider the following ways of defining the
adja\-cency measure.
(N.1) $a_{jk}\,{=}\,I(r_{jk}\,{>}\,r)$ and $s_{jk}\,{=}\,1$. Here the cutoff $r$ is
compu\-ted as 3.09 using the approach described in Section \ref{sec3} with a
$p$-value of $10^{-3}$;
(N.2) $a_{jk}=I(|r_{jk}|>r)$ and $s_{jk}=\sgn(r_{jk})$. Here the
cutoff $r$ is computed as 3.29 using the approach described in Section
\ref{sec3} with a $p$-value of $10^{-3}$;
(N.3) $a_{jk}=\max(0, r_{jk})^{\alpha}$ and $s_{jk}=1$. We set $\alpha
=6$, which satisfies the scale-free topology criteria [\citet{ZhangH05}];
(N.4) $a_{jk}=r_{jk}^{\alpha}$ and $s_{jk}=\sgn(r_{jk})$. We set
$\alpha=6$.

The penalty levels $\lam_1$ and $\lam_2$ are selected using $V$-fold
cross validation. In our numerical study, we set $V=5$. To reduce
computational cost, we search over the discrete grid of $2^{\ldots, -1,
-0.5, 0, 0.5,\ldots}$. For comparison, we also consider the MCP estimate
and the approach proposed in \citet{DayeJeng}; referred to as D--J
hereafter. Both the SLS and MCP involve the regularization parameter
$\gam$. For MCP, \citet{Zhang10} suggested using $\gamma= 2/(1-\max_{j
\neq k }|x_j'x_k|/n)$ for standardized covariates. The average $\gamma
$ value of this choice is 2.69 in his simulation studies. The
simulation studies in Breheny and Huang (\citeyear{BH}) suggest that $\gamma=3$
is a reasonable choice. We have experimented with different $\gamma$
values and reached the same conclusion. Therefore, we set $\gamma=3$.

We set $n=100$ and $p=500$. Among the 500 covariates, there are 100
clusters, each with size 5. We consider two different correlation
structures. (I)~Covariates in different clusters are independent,
whereas covariates $i$ and~$j$ within the same cluster have correlation
coefficients $\rho^{|i-j|}$; and~(II)~covariates $i$ and $j$ have
correlation coefficients $\rho^{|i-j|}$. Under structure I, zero and
nonzero effects are independent, whereas under structure II, they are
correlated. Covariates have marginal normal distributions with mean
zero and variance one. We consider different levels of correlation with
$\rho=0.1, 0.5,\allowbreak 0.9$. Among the 500 covariates, the first 25 (5~clusters)
have nonzero regression coefficients. We consider the
following scenarios for nonzero coefficients: (a)~all the nonzero
coefficients are equal to 0.5; and (b)~the nonzero coefficients are
randomly generated from the uniform distribution on $[0.25, 0.75]$. In~(a),
the Laplacian matrices satisfy the unbiasedness property $L\beta
^o=0$ discussed in Section \ref{sec4}. We have experienced with other levels of
nonzero regression coefficients and reached similar conclusions.

We examine the accuracy of identifying nonzero covariate effects and
the prediction performance. For this purpose, for each simulated
dataset, we simulate an independent testing dataset with sample size
100. We conduct cross validation (for tuning parameter selection) and
estimation\vadjust{\goodbreak} using the training set only. We then make prediction for
subjects in the testing set and compute the PMSE (prediction mean
squared error).

We simulate 500 replicates and present the summary statistics in Table~\ref{table1}.
We can see that the MCP performs satisfactorily when the correlation
is small. However, when the correlation is high, it may miss a
considerable number of true positives and have large prediction errors.
The D--J approach, which can also accommodate the correlation structure,
is able to identify all the true positives. However, it also identifies
a large number of false positives, causing by the over-selection of the
Lasso penalty. The proposed SLS approach outperforms the MCP and D--J
methods in the sense that it has smaller empirical false discovery
rates with comparable false negative rates. It also has significantly
smaller prediction errors.

%
%t1 ###
\begin{sidewaystable}
\tabcolsep=0pt
\textwidth=\textheight
\tablewidth=\textwidth
\caption{Simulation study: median based on 500 replicates. In each
cell, the three numbers are positive findings, true positives and PMSE
$\times100$, respectively}\label{table1}
{\fontsize{8pt}{10pt}\selectfont{\begin{tabular*}{\tablewidth}
{@{\extracolsep{4in minus 4in}}lc@{\hspace*{7pt}}ccd{3.2}@{\hspace*{7pt}}ccd{3.2}@{\hspace*{7pt}}
ccd{3.2}@{\hspace*{7pt}}ccd{3.2}@{\hspace*{7pt}}ccd{3.2}@{\hspace*{7pt}}ccc@{\hspace*{7pt}}ccc@{\hspace*{7pt}}ccc@{\hspace*{7pt}}ccc@{}}
\hline
& & & & & \multicolumn{12}{c}{\textbf{D--J}\hspace*{7pt}}& \multicolumn{12}{c@{}}{\textbf{SLS}}
\\[-4pt]
& & & & & \multicolumn{12}{c}{\hrulefill}& \multicolumn{12}{c@{}}{\hrulefill} \\
\textbf{Coefficient} &$\bolds{\rho}$ & \multicolumn{3}{c}{\textbf{MCP}\hspace*{7pt}} & \multicolumn{3}{c}{\textbf{N.1}\hspace*{7pt}}
& \multicolumn{3}{c}{\textbf{N.2}\hspace*{7pt}} & \multicolumn{3}{c}{\textbf{N.3}\hspace*{7pt}}
& \multicolumn{3}{c}{\textbf{N.4}\hspace*{7pt}} & \multicolumn{3}{c}{\textbf{N.1}\hspace*{7pt}}
& \multicolumn{3}{c}{\textbf{N.2}\hspace*{7pt}} & \multicolumn{3}{c}{\textbf{N.3}\hspace*{7pt}}
& \multicolumn{3}{c@{}}{\textbf{N.4}}\\
%[-4pt]
%& \multicolumn{3}{c}{\hrulefill} & \multicolumn{3}{c}{\hrulefill}
%& \multicolumn{3}{c}{\hrulefill} & \multicolumn{3}{c}{\hrulefill}
%& \multicolumn{3}{c}{\hrulefill} & \multicolumn{3}{c}{\hrulefill}
%& \multicolumn{3}{c@{}}{\hrulefill}\\
\hline\\[-8pt]
& & \multicolumn{27}{c@{}}{Correlation structure I}\\
[4pt]
0.5 & 0.1 & 27 & 25 & 41.33 & 61 & 25 & 125.34 & 53 & 25 & 46.64 & 55 & 25 & 60.14 &
59 & 25 & 51.24 &
27 & 25 & 40.53 & 27 & 25 & 39.84 & 26 & 25 & 41.74 & 27 & 25 & 39.34 \\
& 0.5 & 28 & 25 & 54.10 & 51 & 25 & 66.38 & 67 & 25 & 66.84 & 72 & 25 & 56.22 & 63 & 25
& 53.43 &
27 & 25 & 37.71 & 28 & 25 & 39.18 & 28 & 25 & 33.87 & 27 & 25 & 36.00 \\
& 0.9 & 22 & 15 & 137.52 & 66 & 25 & 55.51 & 55 & 25 & 56.94 & 61 & 25 & 49.22 & 74 & 25
& 51.41 &
29 & 25 & 48.89 & 28 & 25 & 49.96 & 29 & 25 & 45.16 & 27 & 25 & 41.49 \\
$U$[0.25, & 0.1 & 37 & 25 & 52.24 & 72 & 25 & 54.28 & 61 & 25 & 88.00 & 59 & 25
& 70.00 & 78 & 25 & 60.51 &
33 & 25 & 51.80 & 36 & 25 & 52.19 & 30 & 25 & 53.03 & 30 & 25 & 52.22 \\
\quad 0.75] & 0.5 & 29 & 24 & 65.12 & 66 & 25 & 78.76 & 54 & 25 & 72.34 & 63 & 25 & 63.55 & 57 & 25
& 66.33 &
28 & 25 & 42.24 & 28 & 25 & 43.96 & 27 & 24 & 54.72 & 28 & 24 & 58.77 \\
& 0.9 & 17 & 13 & 152.42 & 67 & 25 & 63.43 & 62 & 25 & 57.30 & 50 & 25 & 53.88 & 74 & 25
& 57.98 &
29 & 25 & 47.73 & 29 & 25 & 49.14 & 27 & 25 & 48.49 & 28 & 25 & 50.83 \\[4pt]
& & \multicolumn{27}{c@{}}{Correlation structure II}\\
[4pt]
0.5 & 0.1 & 26 & 25 & 38.22 & 62 & 25 & 121.69 & 58 & 25 & 117.10 & 63 & 25 & 127.34 &
72 & 25 & 122.34 & 27 & 25 & 40.33 & 27 & 25 & 40.65 & 27 & 25 & 41.49 & 27 & 25 & 37.40 \\
& 0.5 & 29 & 25 & 53.01 & 52 & 25 & 55.99 & 49 & 25 & 62.04 & 66 & 25 & 62.70 & 65 & 25
& 64.41 & 27 & 25 & 36.97 & 28 & 25 & 39.47 & 28 & 25 & 38.53 & 27 & 25 & 39.53 \\
& 0.9 & 15 & 13 & 140.69 & 48 & 25 & 55.75 & 34 & 25 & 56.71 & 32 & 25 & 60.27 & 38 & 25
& 59.78 & 29 & 25 & 66.79 & 29 & 25 & 60.52 & 29 & 25 & 57.91 & 30 & 25 & 60.19 \\
$U$[0.25, & 0.1 & 37 & 25 & 54.31 & 77 & 25 & 60.02 & 72 & 25 & 66.14 & 74 & 25
& 78.32 & 66 & 25 & 74.50 & 29 & 25 & 50.05 & 32 & 25 & 51.34 & 37 & 25 & 50.74 & 29 & 25
& 49.47 \\
\quad 0.75] & 0.5 & 27 & 24 & 57.66 & 74 & 25 & 61.71 & 66 & 25 & 67.54 & 75 & 25 & 62.01 & 74 & 25
& 66.91 & 28 & 25 & 44.92 & 28 & 25 & 46.65 & 28 & 25 & 41.35 & 28 & 25 & 41.17 \\
& 0.9 & 14 & 13 & 136.49 & 33 & 25 & 61.50 & 35 & 25 & 55.08 & 34 & 25 & 54.54 & 38 & 25
& 60.67 & 29 & 25 & 56.87 & 29 & 25 & 57.03 & 30 & 25 & 53.28 & 30 & 25 & 56.79 \\
\hline
\end{tabular*}}}
\end{sidewaystable}

%s6.1 ###
\subsection{Application to a microarray study}

In the study reported in \citet{Scheetz}, F1 animals were
intercrossed and 120 twelve-week-old male offspring were selected for
tissue harvesting from the eyes and microarray analysis using the
Affymetric GeneChip Rat Genome 230 2.0 Array. The intensity values were
normalized using the RMA [robust multi-chip averaging, Bolstad et al.
(\citeyear{Boletal03}),
\citet{Iri}] method to obtain summary expression values for each
probe set. Gene expression levels were analyzed on a logarithmic scale.
For the %31,042
probe sets on the array, we first excluded those that were not
expressed in the eye or that lacked sufficient variation. The
definition of expressed was based on the empirical distribution of RMA
normalized values. For a probe set to be considered expressed, the
maximum expression value observed for that probe among the 120~F2 rats
was required to be greater than the 25th percentile of the entire set
of RMA expression values. For a probe to be considered ``sufficiently
variable,'' it had to exhibit at least 2-fold variation in expression
level among the 120 F2 animals.

We are interested in finding the genes whose expression are most
variable and correlated with that of gene TRIM32. This gene was
recently found to cause Bardet--Biedl syndrome [\citet{Chiang}],
which is a genetically heterogeneous disease of multiple organ systems
including the retina. One approach to find the genes related to TRIM32
is to use regression analysis. Since it is expected that the number of
genes associated with gene TRIM32 is small and since we are mainly
interested in genes whose expression values across samples are most
variable, we conduct the following initial screening. We compute the
variances of gene expressions and select the top 1\mbox{,}000. We then
standardize gene expressions to have zero mean and unit variance.

We analyze data using the MCP, D--J, and proposed approach. In cross
validation, we set $V=5$. The numbers of genes identified are MCP: 23,
D--J: 31 (N.1), 41 (N.2), 34 (N.3), 30 (N.4), SLS: 25 (N.1), 26 (N.2),
16~(N.3) and 17 (N.4), respectively. More detailed results are
available from the authors. Different approaches and different ways of
defining the adjacency measure lead to the identification of different
genes. As expected, the SLS identifies shorter lists of genes than the
D--J, which may lead to more parsimonious models and more focused
hypothesis for confirmation. As the proposed approach pays special
attention to the correlation among genes, we also compute the median of
the absolute values of correlations among the identified genes, which
are MCP: 0.171, D--J: 0.201 (N.1), 0.207 (N.2), 0.215 (N.3), 0.206
(N.4), SLS: 0.247 (N.1), 0.208 (N.2), 0.228 (N.3), 0.212 (N.4). The D--J
and SLS, which incorporate correlation in the penalty, identify genes
that are more strongly correlated than the MCP. The SLS identified
genes have slightly higher correlations than those identified by D--J.

Unlike in simulation study, we are not able to evaluate true and false~po\-sitives.
This limitation is shared by most existing studies. We use
the following $V$-fold ($V=5$) cross validation based approach to evaluate
prediction. (a)~Ran\-domly split data into $V$-subsets with equal sizes;
(b)~Remove one subset from data; (c)~Conduct cross validation and
estimation using the rest $V-1$ subsets; (d)~Make prediction for the
one removed subset; (e)~Repeat Steps (b)--(d) over all subsets and
compute the prediction error. The sums of squared prediction errors are
MCP: 1.876; D--J: 1.951 (N.1), 1.694 (N.2), 1.534 (N.3) and 1.528 (N.4);
SLS: 1.842 (N.1), 1.687 (N.2), 1.378 (N.3) and 1.441 (N.4),
respectively. The SLS has smaller cross validated prediction errors,
which may indirectly suggest better selection properties.

%s7 ###
\vspace*{3pt}\section{Discussion}\label{sec7}
In this article, we propose the SLS method for variable selection and
estimation in high-dimensional data analysis. The most important
feature of the SLS is that it explicitly incorporates the graph/network
structure in predictors into the variable selection procedure through
the Laplacian quadratic. It provides a systematic framework for
connecting penalized methods for consistent variable selection and
those for network and correlation analysis. As can be seen from the
methodological development, the application of the SLS variable
selection is relatively independent of the graph/network construction.
Thus, although graph/network construction is of significant importance,
it is not the focus of this study and not thoroughly pursued.

An important feature of the SLS method is that it incorporates the
correlation patterns of the predictors into variable selection through
the Laplacian quadratic. We have considered two simple approaches for
determining the Laplacian based on dissimilarity and similarity
measures. Our simulation studies demonstrate that incorporating
correlation patterns improves selection results and prediction
performance. Our theoretical results on the selection properties of the
SLS are applicable to a general class of Laplacians and do not require
the underlying graph for the predictors to be correctly specified.

We provide sufficient conditions under which the SLS estimator
possesses an oracle property, meaning that it is sign consistent and
equal to the oracle Laplacian shrinkage estimator with high
probability. We also study the grouping properties of the SLS
estimator. Our results show that the SLS is adaptive to the sparseness
of the original $p$-dimensional model with $p \gg n$ and the denseness
of the underlying $d^o$-dimensional model, where $d^o < n $ is the
number of nonzero coefficients. The asymptotic rates of the penalty
parameters are derived. However, as in many recent studies, it is not
clear whether the penalty parameters selected using cross validation or
other procedures can match the asymptotic rate. This is an important
and challenging problem that requires further investigation, but is
beyond the scope of the current paper. Our numerical study shows a
satisfactory finite-sample performance of the SLS. Particularly, we
note that the cross validation selected tuning parameters seem
sufficient for our simulated data. We are only able to experiment with
four different adjacency measures. It is not our intention to draw
conclusions on different ways of defining adjacency. More adjacency
measures are hence not explored.

We have focused on the linear regression model in this article.
However, the SLS method can be applied to general linear regression
models. Specifically, for general linear models, the SLS criterion can
be formulated as
\begin{eqnarray*}
&&\frac{1}{n} \sum_{i=1}^n \ell\biggl(y_i, b_0+\sum_{j}x_{ij}b_j\biggr) +
\sum_{j=1}^p\rho(|b_j|; \lam_1, \gamma) \\
&&\qquad{} + \frac{1}{2}\lam_2 \sum
_{1\le j < k \le p} \vert a_{jk}\vert(b_j-s_{jk}b_k)^2,
\end{eqnarray*}
where $\ell$ is a given loss function. For instance, for generalized
linear models such as logistic regression, we can take $\ell$ to be
the negative log-likelihood function. For Cox regression, we can use
the negative partial likelihood as the loss function. Computationally,
for loss functions other than least squares, the coordinate descent
algorithm can be applied iteratively to quadratic approximations to the
loss function. However, further work is needed to study theoretical
properties of the SLS estimators for general linear models.

There is a large literature on the analysis of network data
and much work has also been done
on estimating sparse covariance matrices in high-dimensional settings.
See, for example, \citet{ZhangH05}, \citet{Chungb},
\citet{MeinshausenB06}, \citet{YuanLin},
\citet{FHT}, \citet{FFW}, among others.
It would be useful to study ways to incorporate these methods and
results into the proposed SLS approach.
In some problems such as genomic data analysis, partial external
information may also be available on the graphical structure of some
genes used as predictors in the model. It would be interesting to
consider approaches for combining external information on the graphical
structure with existing data in constructing the Laplacian quadratic penalty.

\begin{appendix}\label{app}
\vspace*{3pt}\section*{Appendix}\vspace*{3pt}
In this Appendix, we give proofs of Theorems \ref{ThmA} and \ref{ThmB}
and Propositions~\ref{PropA} and~\ref{PropB}.
\begin{pf*}{Proof of Theorem \ref{ThmA}}
Since\vspace*{1pt} $c_{\min}(\lam_2) > 1/\gamma$, the criterion (\ref{SLSdefA})
is strictly convex
and its minimizer is unique.
Let $\tX=\tX(\lam_2)=\sqrt{n}(\Sigma+\lam_2L)^{1/2}$, $\tby=\tby
(\lam_2)=\tX^{-1}X'\by$ and
\[
\tM(\bb;\lam,\gamma) = (2n)^{-1}\|\tby-\tX\bb\|^2+\sum
_{j=1}^p\rho(|b_j|;\lam_1,\gamma).
\]
Since\vspace*{1pt} $\tX'(\tX/n,\tby)=(\Sigma+\lam_2L,X'\by)$, $M(\bb;\lam
,\gamma)-\tM(\bb;\lam,\gamma)
=(\|\by\|^2-\|\tby\|^2)/\allowbreak(2n)$ does not depend on $\bb$. Thus,
$\hbbeta$ is the minimizer of
$\tM(\bb;\lam,\gamma)$.

Since $|\hbeta{}^o_j| \ge\gamma\lambda_1$ gives $\rho'(|\hbeta
{}^o_j|;\lambda_1)=0$,
the KKT conditions hold for $\tM(\bb;\allowbreak\lam,\gamma)$ at
$\hbbeta(\lam)=\hbbeta{}^o(\lam)$ in the intersection of the events
%
%e7.1 ###
\begin{equation}
\Omega_1 =\{\|\tX'_{\cO^c}(\tby-\tX\hbbeta{}^o)/n\|_\infty\le\lam
_1 \},\quad
\Omega_2 = \Bigl\{ \min_{j \in\cO} \sgn(\beta_j^*)\hbeta{}^o_j \ge
\gamma\lam_1\Bigr\}.\hspace*{-40pt}
\label{KKTb}
\end{equation}
Let $\tbveps^* = \tby- \tX\bbeta^* = \tbveps+ E\tbveps^*$ with
$\tbveps= \tby- E\tby$. Since $\tX'\tby=X'\by$ and both~$\bbeta^o$
and $\bbeta^*$ are supported in $\cO$,
%
%e7.2 ###
\begin{eqnarray}\label{bias-1} \tX_B'E\tbveps^*/n
&=& X_B'X\bbeta^o/n-\tX_B'\tX\bbeta^*/n \nonumber\\
&=&
\Sigma_{B,\mathcal{O}}\bbeta^o_\mathcal{O} -
\Sigma_{B,\mathcal{O}}(\lam_2)\Sigma_{\mathcal{O}}^{-1}(\lam
_2)\Sigma_\mathcal{O}\bbeta^o_\mathcal{O}
\\
&=&
\lam_2\{\Sigma_{B,\mathcal{O}}(\lam_2)\Sigma_\mathcal
{O}^{-1}(\lam_2)L_\mathcal{O}-L_{B,\mathcal{O}}\}
\bbeta^o_\mathcal{O},\nonumber
\end{eqnarray}
which\vspace*{1pt} describes the effect of the bias of
$\hbbeta{}^o$ on the gradient in the linear model $\tby=
\tX\bbeta^*+\tbveps^*$. Since $\tX_\mathcal{O}'E\tbveps^*/n=0$,
we have
$\|\tX'E\tbveps^*/n\|_\infty=\lam_2C_2$.

Since $\tX'\tbveps=\tX'\tby-E\tX'\tby= X'\by-EX'\by=X'\bveps$,
(\ref{bias-1}) gives
%
%e7.3 ###
\begin{equation}\label{Omega1}
\Omega_1\subseteq\{\|X_{\mathcal{O}^c}'\bveps/n\|_\infty<
\lam_1-\lam_2C_2\}.
\end{equation}
Since $\bbeta^*=E\hbbeta{}^o$, $\hbbeta{}^o_\mathcal{O}=\Sigma
_\mathcal{O}^{-1}(\lam_2)X_\mathcal{O}'\by/n$
can be written as $\bbeta^*_\mathcal{O}+((v_j/n)^{1/2}\bu_j'\bveps$,
$j\,{\in}\,\mathcal{O})'$,
where $\|\bu_j\|\,{=}\,1$ and $\{v_j, j\,{\in}\,\mathcal{O}\}$ are the diagonal
elements of
$\Sigma_\mathcal{O}^{-1}(\lam_2)\,{\times}\,\allowbreak\Sigma_\mathcal{O}\{\Sigma
_\mathcal{O}^{-1}(\lam_2)\}$. Thus,
%
%e7.4 ###
\begin{equation}\label{Omega2}
\Omega_2^c \subseteq\bigcup_{j\in\mathcal{O}}\bigl\{\sgn(\beta
^*_j)\bu_j'\bveps\ge
(n/v_j)^{1/2}|\beta^*_j|\ge\sigma\sqrt{2\log(|\mathcal{O}|/\eps
)}\bigr\}.
\end{equation}
Since $\lam_1\ge\lam_2C_2
+ \sigma\sqrt{2\log(p/\eps)}\max_{j\le p}\|\bx_j\|/n$, the
sub-Gaussian Condition~\ref{condA} yields
\begin{eqnarray*}
1 - P\{\Omega_1\cap\Omega_2\}
&\le& P\Bigl\{\|X_{\cO^c}'\bveps/n\|_\infty> \sigma\sqrt{2\log
\bigl((p-|\cO|)/\eps\bigr)}\max_{j\le p}\|\bx_j\|/n\Bigr\}
\\
&&{} + \sum_{j\in\mathcal{O}}P\bigl\{\sgn(\beta^*_j)\bu
_j'\bveps\ge\sigma\sqrt{2\log(|\mathcal{O}|/\eps)}\bigr\}
\\
&\le& 2|\cO^c|\eps/(p-|\cO|) + |\mathcal{O}|\eps/|\mathcal
{O}|=3\eps.
\end{eqnarray*}
The proof of (\ref{ThmA-1}) is complete, since $\hbeta{}^o_j\neq0$ for
all $j\in\mathcal{O}$ in $\Omega_2$.

For the proof of (\ref{ThmA-2}), we have $\|\bbeta^*_\mathcal{O} -
\bbeta^o_\mathcal{O}\|_\infty=\lam_2C_1$
due to
%
%e7.5 ###
\begin{equation}\label{bias-2}
\bbeta^*_\mathcal{O} - \bbeta^o_\mathcal{O}
= \Sigma_{\mathcal{O}}^{-1}(\lam_2)\Sigma_\mathcal{O}\bbeta
^o_\mathcal{O}- \bbeta^o_\mathcal{O}
= - \lam_2\Sigma_{\mathcal{O}}^{-1}(\lam_2)L_\mathcal{O}\bbeta
^o_\mathcal{O}.
\end{equation}
It follows that the condition on $\beta_*$ implies Condition \ref{condB}(iii) with
$\sgn(\bbeta^*_{\cO})=\sgn(\bbeta^o_{\cO})=\sgn(\hbbeta{}^o_{\cO
})$ in $\Omega_2$.
\end{pf*}
\begin{pf*}{Proof of Theorem \ref{ThmB}} For $m \ge1$ and vectors $\bu
$ in the range of $\tX$, define
%
%e7.6 ###
\begin{eqnarray}
\label{tzeta}
&&\tzeta(\bv;m,\cO,\lambda_2) \nonumber\\[-8pt]\\[-8pt]
&&\qquad= \max\biggl\{\frac{\|(\tP_B-\tP
_{{\cO}})\bv\|_2}{(m n)^{1/2}}\dvtx
\cO\subseteq B\subseteq\{1,\ldots. p\}, |B| = m+|\cO|
\biggr\},\nonumber
\end{eqnarray}
where $\tP_B = \tX_B (\tX_B'\tX_B)^{-1} \tX_B'$.
Here $\tzeta$ depends on $\lam_2$ through $\tP$.
Since $\hbbeta(\lam)$ is the MC+ estimator based on data $(\tX,\tby
)$ at penalty level $\lam_1$
and (\ref{SRCa}) holds for $\Sigma(\lam_2)=\tX'\tX/n$, the proof
of Theorem 5\vspace*{1pt} in \citet{Zhang10} gives
$\hbbeta(\lam)=\hbbeta{}^o(\lam)$ in the event $\Omega=\bigcap
_{j=1}^3\Omega_j$, where
$\Omega_1= \{\|\tX_{\mathcal{O}^c}'(\tby-\tX\hbbeta{}^o)/n\|
_\infty\le\lam_1\}$ is as in (\ref{KKTb})
and
\begin{eqnarray*}
\Omega_2 &=& \Bigl\{\min_{j\in\mathcal{O}}\sgn(\beta_j^*)\hbeta{}^o_j
> \gamma\bigl(2\sqrt{c^*}\lam_1\bigr)\Bigr\},\\
\Omega_3 &=& \{\zeta(\tby-\tX\bbeta^*;d^*-|\mathcal
{O}|,\mathcal{O},\lam_2)\le\lam_1\}.
\end{eqnarray*}
Note that $(\lam_{1,\eps},\lam_{2,\eps},\lam_{3,\eps},\alpha)$
in \citet{Zhang10}
is identified with $(\lam_1,2\sqrt{c^*}\lam_1,\lam_1$, $1/2)$ here.

Let\vspace*{1pt} $\tbveps^* = \tby- \tX\bbeta^* = \tbveps+ E\tbveps^*$ with
$\tbveps= \tby- E\tby$.
Since $\tX'\tby=X'\by$, (\ref{bias-1}) still holds with
$\|\tX'E\tbveps^*/n\|_\infty=\lam_2C_2$.
Since $\tX'\tbveps=X'\by-EX'\by=X'\bveps$, (\ref{bias-1}) still
gives (\ref{Omega1}).
A slight modification of the argument for (\ref{Omega2}) yields
%
%e7.7 ###
\begin{eqnarray}\label{Omega2a}
\Omega_2^c &\subseteq&\bigcup_{j\in\mathcal{O}}\bigl\{\sgn(\beta
^*_j)\bu_j'\bveps\ge
(n/v_j)^{1/2}\bigl(|\beta^*_j| - \gamma\bigl(2\sqrt{c^*}\lam_1\bigr)\bigr)\nonumber\\[-8pt]\\[-8pt]
&&\hspace*{137.6pt}\ge
\sigma\sqrt{2\log(|\mathcal{O}|/\eps)}\bigr\}.\nonumber
\end{eqnarray}
For $|B|\le d^*$, we have $\|\tP_BE\tbveps^*\|/\sqrt{n}=\|\Sigma
_B^{-1/2}(\lam_2)\tX_B'E\tbveps^*\|/ n
\le\break\|\tX_B'E\tbveps^*/ n\|_\infty\sqrt{|B|/c_*(\lam_2)}$ and
$\|\tP_B\tbveps\|/\sqrt{n}=\|\Sigma_B^{-1/2}(\lam_2)\tX_B'\tbveps
\|/n\le\break\|X_B'\bveps/n\|_\infty\times\sqrt{|B|/c_*(\lam_2)}$. Thus, by (\ref{tzeta})
\begin{eqnarray*}
\zeta(\tby-\tX\bbeta^*;d^*-|\mathcal{O}|,\mathcal{O},\lam_2)
&=& \zeta(\tbveps+E\tbveps^*;d^*-|\mathcal{O}|,\mathcal{O},\lam_2)\\
&\le&\frac{(\|X'\bveps/n\|_\infty+\lam_2C_2)\sqrt{d^*}}{\sqrt
{(d^*-|\mathcal{O}|)c_*(\lam_2)}}.
\end{eqnarray*}
Since $|\mathcal{O}|\le d^*/(K_*+1)$, this gives
%
%e7.8 ###
\begin{equation}\label{Omega3}
\Omega_3\subseteq\bigl\{\|X'\bveps/n\|_\infty< \sqrt{c_*(\lam
_2)K_*/(K_*+1)}\lam_1-\lam_2C_2\bigr\}.
\end{equation}
Since $\max\{1,\sqrt{c_*(\lam_2)K_*/(K_*+1)}\}\lam_1\ge\lam_2C_2
+ \sigma\sqrt{2\log(p/\eps)}\max_{j\le p}\|\bx_j\|/n$,
(\ref{Omega1}), (\ref{Omega2a}), (\ref{Omega3}) and Condition \ref{condA}
imply
\begin{eqnarray*}
&& 1 - P\{\Omega_1\cap\Omega_3\} + P\{\Omega_2^c
\}
\\
&&\qquad\le P\Bigl\{\|X'\bveps/n\|_\infty> \sigma\sqrt{2\log(p/\eps
)}\max_{j\le p}\|\bx_j\|\big/n\Bigr\}\\
&&\qquad\quad{}+ \sum_{j\in\mathcal{O}}P\bigl\{\sgn(\beta^*_j)\bu_j'\bveps\ge
\sigma\sqrt{2\log(|\mathcal{O}|/\eps)}\bigr\}
\\
&&\qquad\le 2p(\eps/p) + |\mathcal{O}|\eps/|\mathcal{O}|=3\eps.
\end{eqnarray*}
The proof of (\ref{ThmB-1}) is complete, since $\hbeta{}^o_j\neq0$ for
all $j\in\mathcal{O}$ in $\Omega_2$.
We omit the proof of (\ref{ThmB-2}) since it is identical to that of
(\ref{ThmA-2}).
\end{pf*}
\begin{pf*}{Proof of Proposition \ref{PropA}}
The $\tbbeta$ satisfies
%
%e7.9 ###
\begin{equation}
\label{tsa}
-\frac{1}{n}\bx_j'(\by-X\tbbeta)+\lam_2 (d_j\tbeta_j-\ba
_j'\tbbeta)=0, \qquad 1 \le j \le q.
\end{equation}
Therefore, by Cauchy--Schwarz and using $\|\bx_j\|^2=n$, we have
\[
\lam_2 \max_{1 \le j \le q} |d_j\tbeta_j-\ba_j'\tbbeta|
\le\frac{1}{n}\max_{1\le j \le q}|\bx_j'(\by-X\tbbeta)|
\le\frac{1}{\sqrt{n}}\|\tbr\|.
\]
Now because $G(\tbbeta;\lam_2) \le G(\bzero;\lam_2)$,
we have $\|\tbr\| \le\|\by\|$. This proves part (i).

For part (ii), note that we have
\[
\lam_2\bigl(d_j\tbeta_j-\ba_j'\tbbeta- (d_k\tbeta_k-\ba_k'\tbbeta)\bigr)
=\frac{1}{n}(\bx_j-\bx_k)'\tbr.
\]
Thus
\[
\lam_2|d_j\tbeta_j-\ba_j'\tbbeta- (d_k\tbeta_k-\ba_k'\tbbeta)|
\le\frac{1}{n} \|\bx_j-\bx_k\| \|\tbr\|.
\]
%
%Since $G(\tbbeta;\lam_2) \le G(0;\lam_2)$, we have $\|\tbr\|\le\|\by
Part (ii) follows.
\end{pf*}
\begin{pf*}{Proof of Proposition \ref{PropB}}
The $\tbbeta$ must satisfy
%
%e7.10 ###
\begin{equation}
\label{Neq1}\qquad
-\frac{1}{n}\bx_j'(\by-X\tbbeta)+\lam_2
(\tbeta_j-v_g^{-1}\bone_g'\tbbeta_g)=0, \qquad j \in V_g,  1\le g \le J.
\end{equation}
Taking the difference between the $j$th and $k$th equations in (\ref
{Neq1}) for \mbox{$j, k \in V_g$}, we get
\[
\lam_2(\tbeta_j-\tbeta_k)=\frac{1}{n}(\bx_j-\bx_k)'(\by-X\tbbeta
),  \qquad j, k \in V_g.
\]
Therefore,
\[
\lam_2\vert\tbeta_j-\tbeta_k\vert \le\frac{1}{n}\|\bx_j-\bx_k\|
\cdot\|\by-X\tbbeta\|, \qquad  j, k \in V_g.
\]
Part (i) follows from this inequality.\vspace*{1pt}
%Thus for any $\lam> 0$, when $\|\bx_j-\bx_k\| \rightarrow0$, we have
%$\abs{\tbeta_j-\tbeta_k} \rightarrow0$.

Define\vspace*{1pt} $\bar{\beta}_g=v_g^{-1}\bone_g'\tbbeta_g$. This is the
average of the elements in $\tbbeta_g$.
For any $j \in V_g$ and $k \in V_h, g\neq h$, we have
\[
\lam_2\bigl(\tbeta_j-\bar{\beta}_g-(\tbeta_k-\bar{\beta}_h)\bigr)=\frac
{1}{n}(\bx_j-\bx_k)'(\by-X\tbbeta), \qquad  j \in V_g, k \in V_h.
\]
Thus, part (ii) follows. This completes the proof of Proposition \ref
{PropB}.
\end{pf*}
\end{appendix}

\section*{Acknowledgments}
We wish to thank two anonymous referees, the Associate Editor and
Editor for their helpful comments which led to considerable
improvements in the presentation of the paper.

%suskaldyti doi

% imsref loaded by lrinkeviciute, 2011-06-14 09:38:10
% imsref loaded by lrinkeviciute, 2011-06-14 09:44:29
%
% imsref loaded by lrinkeviciute, 2011-06-14 13:48:50

%
\printaddresses

\end{document}